\documentclass[12pt,reqno]{amsart}

\usepackage{amssymb}

\usepackage{eucal}

\input{diagrams}

\font\sss=cmss8

\def\sD{\mbox{\sf D}}

\def\ast{{\textstyle *}}

\def\D{\sD}

\def\Dfd{\sD^{\operatorname{fd}}}
\def\Did{\sD^{\operatorname{id}}}

\def\Dsmall{\mbox{\sss D}}

\def\Ext{\operatorname{Ext}}

\def\fd{\operatorname{fd}}

\def\FPD{\operatorname{FPD}}

\def\H{\operatorname{H}}

\def\Hom{\operatorname{Hom}}
\def\id{\operatorname{id}}

\def\inf{\operatorname{inf}}

\def\Ksmall{\mbox{\sss K}}

\def\LTensor{\stackrel{\operatorname{L}}{\otimes}}

\def\opp{\operatorname{op}}
\def\pd{\operatorname{pd}}

\def\RHom{\operatorname{RHom}}

\def\sup{\operatorname{sup}}

\numberwithin{equation}{part}

\newtheorem{Lemma}{Lemma}[section]

\newtheorem{Theorem}[Lemma]{Theorem}

\newtheorem{Corollary}[Lemma]{Corollary}

\theoremstyle{definition}
\newtheorem{Definition}[Lemma]{Definition}
\newtheorem{Setup}[Lemma]{Setup}

\begin{document}

\title[Flat and projective dimension]
{Finite flat and projective dimension}

\author{Peter J\o rgensen}
\address{Department of Pure Mathematics, University of Leeds,
Leeds LS2 9JT, United Kingdom}
\email{popjoerg@maths.leeds.ac.uk, www.maths.leeds.ac.uk/\~{ }popjoerg}


\keywords{Finitistic projective dimension, flat dimension, dualizing
complex, complete semi-local PI algebra, noetherian connected
filtration}

\subjclass[2000]{16E10}



\maketitle

It is surprising but true that if $A$ is a nice enough ring,
then any $A$-left-module of finite flat dimension also has
finite projective dimension.

A result to this effect was proved in \cite[prop.\ 6]{Jensen} under
the condition that the number
\[
  \FPD(A) =
  \sup 
  \biggl\{\, 
  \pd M \,\,\bigg|
  \begin{array}{l}
    \mbox{ $M$ is an $A$-left-module } \\
    \mbox{ with $\pd M < \infty$ }
  \end{array}
  \!\!\biggr\},
\]
known as the left big finitistic projective dimension of $A$, is
finite.  Here $\pd M$ denotes the projective dimension of $M$.

Unfortunately, this number is not known to be finite even if $A$ is a
finite dimensional algebra over a field, where, indeed, its finiteness
is a celebrated conjecture.  On the other hand, for such an algebra,
finite flat certainly implies finite projective dimension, simply
because each flat module is projective.

So it seems that there might be results based on other conditions than
finiteness of $\FPD(A)$.  The germ of such a result is in
\cite[cor.\ 3.4]{Foxby} which shows that finite flat implies finite
projective dimension for a ring which is a homomorphic image of a
noetherian commutative Gorenstein ring with finite Krull dimension.

Now, for such a ring, $\FPD(A)$ is in fact finite, and so formally,
\cite[cor.\ 3.4]{Foxby} follows from \cite[prop.\ 6]{Jensen}.
However, the method of proof of \cite[cor.\ 3.4]{Foxby} lends itself
to generalization, and in this note I will use it to show that finite
flat implies finite projective dimension for any right-noetherian
algebra which admits a dualizing complex.  

This includes finite dimensional algebras, but also more interesting
cases such as noetherian complete semi-local PI algebras, and filtered
algebras whose associated graded algebras are connected and
noetherian, and either PI, graded FBN, or with enough normal elements.

\begin{Setup}
$k$ is a field, and $A$ is a right-noetherian $k$-algebra for which
there exists a left-noetherian $k$-algebra $B$ and a
dualizing complex ${}_{B}D_A$.
\end{Setup}

To make sense of this setup, here is the definition of dualizing
complexes, reproduced from \cite[def.\ 1.1]{YekutieliZhang}.  Note
that the definition uses derived categories; for notation relating to
these, see for instance \cite[sec.\ 1]{WuZhangDualizing} or
\cite[sec.\ 1]{YekutieliZhang}.

\begin{Definition}
The complex ${}_{B}D_A$ in $\D(B \otimes_k A^{\opp})$ is called a
dualizing complex if it satisfies the following.
\begin{enumerate}

  \item  The cohomology of $D$ is bounded and finitely generated
         both over $B$ and over $A^{\opp}$.

  \item  The injective dimensions $\id_B D$ and $\id_{A^{\opp}} D$ are
         finite. 

  \item  The canonical morphisms $A \longrightarrow \RHom_B(D,D)$ in
         $\D(A^e)$ and $B \longrightarrow \RHom_{A^{\opp}}(D,D)$ in
         $\D(B^e)$ are isomorphisms.

\end{enumerate}
\end{Definition}

The dualizing complex $D$ is a sort of co-tilting object, so it is not
surprising that the functors $D \LTensor -$ and $\RHom(D,-)$ induce
certain equivalences of categories.  These were noted in \cite[sec.\
2]{AvrFoxPLMS} in the commutative case and in \cite[sec.\ 2]{Mori} in
the non-commutative graded case.  The methods carry over to the
present setup with the following result: Let $\Dfd(A)$ denote the
full subcategory of $\D(A)$ consisting of complexes which are
isomorphic to a bounded complex of flat modules, and let $\Did(B)$
denote the full subcategory of $\D(B)$ consisting of complexes which
are isomorphic to a bounded complex of injective modules.  Then there
are quasi-inverse equivalences of categories
\begin{equation}
\label{equ:equivalences}
  \begin{diagram}[labelstyle=\scriptstyle,width=12ex]
    \Dfd(A) & \pile{ \rTo^{D \LTensor_A -} 
                     \\ \lTo_{\RHom_B(D,-)} } 
                & \Did(B). \\
  \end{diagram}
\end{equation}

Using this, I can prove the main result that finite flat implies
finite projective dimension.

\begin{Theorem}
Let $M$ be an $A$-left-module.  Then $\fd M < \infty$ implies $\pd M <
\infty$. 
\end{Theorem}

\begin{proof}
The case $M = 0$ is trivial, so let me assume $M \not\cong 0$.  I will
then prove more than claimed, namely, $\fd M < \infty$ implies
\begin{equation}
\label{equ:inequality}
  \pd_A M
  \leq  -\inf \{\, i \,|\, \H^i(D \LTensor_A M) \not\cong 0 \,\} + \id_B D.
\end{equation}
This implies the theorem because both terms on the right hand side are
$< \infty$, the first one since $D$ has bounded cohomology
while $M$ has finite flat dimension, and the second one by assumption
on $D$.

The right hand side of \eqref{equ:inequality} does not change if $D$
is replaced by some suspension $\Sigma^j D$, so I can suppose
\[
  \id_B D = 0.
\]
If I now set
\[
  n  =  -\inf \{\, i \,|\, \H^i(D \LTensor_A M) \not\cong 0 \,\}
\]
then the inequality \eqref{equ:inequality} amounts to
\[
  \pd_A M \leq n.
\]

To prove this, let me start by showing $n \geq 0$.  The condition
$\id_B D = 0$ implies that the cohomology of $D$ is concentrated in
cohomological degrees $\leq 0$, and since $M$ is a complex
concentrated in degree $0$ and so has cohomology concentrated in
degree $0$, it follows that the cohomology $\H(D \LTensor_A M)$ is
concentrated in degrees $\leq 0$.  Provided 
$\H(D \LTensor_A M) \not\cong 0$ holds,
this proves $\inf \{\, i \,|\, \H^i(D \LTensor_A M) \not\cong 0
\,\} \leq 0$ and hence $n \geq 0$.  To see 
$\H(D \LTensor_A M) \not\cong 0$,
observe that as $M$ has finite flat dimension, the
equivalences \eqref{equ:equivalences} give the isomorphism $M \cong
\RHom_B(D,D \LTensor_A M)$.  So $M \not\cong 0$ implies $D \LTensor_A
M \not\cong 0$, hence $\H(D \LTensor_A M) \not\cong 0$.

Now consider a projective resolution $P$ of $M$.  There is a short
exact sequence
\[
  0 \rightarrow \Omega^{n+1}M 
    \longrightarrow P_n 
    \longrightarrow \Omega^n M 
    \rightarrow 0
\]
involving two of the syzygies of $M$, as defined by means of $P$.  I
shall prove that this sequence splits, whence $\Omega^n M$ is
projective so $\pd_A M \leq n$ as desired.

To see that the sequence splits, I will in fact prove
\[
  \Ext_A^1(\Omega^n M,\Omega^{n+1}M) = 0.  
\]
For this, consider the following chain of isomorphisms,
\begin{eqnarray*}
  \Ext_A^1(\Omega^n M,\Omega^{n+1}M)
  & \cong & \Ext_A^{n+1}(M,\Omega^{n+1}M) \\
  & \cong & \Hom_{\Dsmall(A)}(M,\Sigma^{n+1}\Omega^{n+1}M) \\
  & \stackrel{\rm (a)}{\cong} 
    & \Hom_{\Dsmall(B)}(D \LTensor_A M,
                        D \LTensor_A (\Sigma^{n+1}\Omega^{n+1}M)) \\
  & = & (*),
\end{eqnarray*}
where $\Sigma^{n+1}$ denotes $(n+1)$'st suspension, and where (a) is
by the equivalences in equation \eqref{equ:equivalences} which apply
because both $M$ and $\Omega^{n+1}M$ have finite flat dimension.

Let $T$ be a truncation of $D \LTensor_A M$ at cohomological degree
\[
  -n = \inf \{\, i \,|\, \H^i(D \LTensor_A M) \not\cong 0 \,\},
\]
so $T$ is quasi-isomorphic to $D \LTensor_A M$ and concentrated in
degrees $\geq -n$.  Then
\[
  (*) \cong \Hom_{\Dsmall(B)}(T,
                              D \LTensor_A (\Sigma^{n+1}\Omega^{n+1}M))
  = (**).
\]

To continue the computation, let ${}_{B}I_A$ be a left-bounded
injective resolution of ${}_{B}D_A$ over $B \otimes_k A^{\opp}$.
Forgetting the $A$-structure, ${}_{B}I$ is an injective resolution of
${}_{B}D$ over $B$.  The condition $\id_B D = 0$ implies that
truncating ${}_{B}I$ at cohomological degree $0$ gives a
quasi-isomorphic complex which is concentrated in degrees $\leq 0$ and
consists of injective $B$-left-modules.  So truncating ${}_{B}I_A$ at
degree $0$ gives a quasi-isomorphic complex ${}_{B}J_A$ which is
concentrated in degrees $\leq 0$ and consists of
$B$-left-$A$-right-modules which are injective when viewed as
$B$-left-modules.  As $I$ is left-bounded, the truncation $J$ is
bounded.

Also, let $F$ be a bounded flat resolution of
$\Sigma^{n+1}\Omega^{n+1}M$.  As the complex
$\Sigma^{n+1}\Omega^{n+1}M$ is just the $(n+1)$'st suspension of a
module, $F$ can be taken to be concentrated in cohomological degrees
$\leq -n-1$, and I have
\[
  D \LTensor_A (\Sigma^{n+1}\Omega^{n+1}M)
  \cong  D \otimes_A F
  \simeq J \otimes_A F.
\]

As $J$ is a bounded complex of $B$-left-$A$-right-modules which are
injective when viewed as $B$-left-modules while $F$ is a bounded
complex of flat $A$-left-modules, $J \otimes_A F$ is a bounded complex
of injective $B$-left-modules.  Hence $J \otimes_A F$ is an injective
resolution which can be used to compute homomorphism groups in
$\D(B)$, so 
\[
  (**) \cong \Hom_{\Ksmall(B)}(T,J \otimes_A F) = (\ast\ast\ast).
\]

But $J$ is concentrated in cohomological degrees $\leq 0$ and $F$ is
concentrated in degrees $\leq -n-1$, so $J \otimes_A F$ is
concentrated in degrees $\leq -n-1$.  On the other hand, $T$ is
concentrated in degrees $\geq -n$.  This implies
\[
  (\ast\ast\ast) = 0
\]
as desired.
\end{proof}

\begin{Corollary}
Let $R$ be a noetherian $k$-algebra and suppose that one of the
following holds.
\begin{enumerate}

  \item  $R$ is a complete semi-local PI algebra.

  \item  $R$ has a filtration $F$ so that the associated graded
         algebra $\operatorname{gr}^F\!R$ is connected and noetherian,
         and either PI, graded FBN, or with enough normal elements.

\end{enumerate}
Let $M$ be an $R$-left-module.  Then $\fd M < \infty$ implies $\pd M <
\infty$.
\end{Corollary}

\begin{proof}
The algebra $R$ can be used as $A$ in the theorem because $B$ and $D$
exist.  In case (i) this is by \cite[cor.\ 0.2]{WuZhangDualizing}, and
in case (ii) by \cite[cor.\ 6.9]{YekutieliZhang}.
\end{proof}

\bigskip

\noindent
{\bf Acknowledgement.}  The diagrams were typeset with Paul Taylor's
{\tt diagrams.tex}.

\end{document}